\documentclass[final,1p,times]{elsarticle}

\usepackage{amssymb}
\usepackage{amsthm}
\usepackage{amscd}
\usepackage{bm}
\usepackage{amsmath}
\usepackage{amsfonts}
\usepackage{graphicx}
\newtheorem{theorem}{Theorem}[section]
\newtheorem{remark}[theorem]{Remark}

\newtheorem{lemma}[theorem]{Lemma}

\usepackage{mathrsfs}
\usepackage{titletoc}
\usepackage{color}

\newcommand{\p}{\partial}
\newcommand{\f}{\frac}
\newcommand{\g}{\gamma}

\newcommand{\be}{\begin{equation}}
\newcommand{\ee}{\end{equation}}
\newcommand{\bea}{\begin{eqnarray}}
\newcommand{\eea}{\end{eqnarray}}
\newcommand{\bna}{\begin{eqnarray*}}
\newcommand{\ena}{\end{eqnarray*}}

\renewcommand{\le}{\left}

\usepackage{geometry}
\journal{***}

\begin{document}

\begin{frontmatter}

\title{Existence and convergence of solutions to $p$-Laplace equations on locally finite graphs}

\author[qnu]{Mengqiu Shao}
\ead{mqshaomath@qfnu.edu.cn}
\address[qnu]{School of Mathematical Sciences, Qufu Normal University, Shandong, 273165, China}

\author[ruc]{Yunyan Yang}
\ead{yunyanyang@ruc.edu.cn}
\address[ruc]{School of Mathematics, Renmin University of China, Beijing, 100872, China}
\author[bnu]{Liang Zhao\corref{Zhao}}
\ead{liangzhao@bnu.edu.cn}
\address[bnu]{School of Mathematical Sciences, Key Laboratory of Mathematics and Complex Systems of MOE,\\
Beijing Normal University, Beijing, 100875, China\\ \vspace{.5cm} {\small\it Dedicated to Professor Jiayu Li on the occasion of his 60th birthday}}

\cortext[Zhao]{Corresponding author.}

\begin{abstract}
We are mainly concerned with the nonlinear $p$-Laplace equation
\begin{equation*}
	-\Delta_pu+\rho|u|^{p-2}u=\psi(x,u)
\end{equation*}
on a locally finite graph $G=(V,E)$, where $p$ belongs to $(1, +\infty)$.
We obtain existence of positive solutions and positive ground state solutions by using the mountain-pass theorem and the Nehari manifold respectively. Moreover, we also analyze the asymptotic behavior for a sequence of positive ground state solutions.
Compared with all the existing relevant works, our results have made essential improvements in at least three aspects: $(i)$ $p$ can take any value in $(1, +\infty)$; $(ii)$ the conditions on the graph $G$ and the potential $\rho$ are relaxed; $(iii)$ for the existence of positive solutions, the growth condition in previous works on the nonlinear term $\psi(x,s)$ as $s\rightarrow +\infty$ is removed.
\end{abstract}


\end{frontmatter}

\section{Introduction}

The variational method for nonlinear PDEs on graphs was pioneered by Grigor'yan, Lin and Yang in their series of works \cite{Gri1,Gri2,Gri3}. In recent years, it has become a very attractive topic and there have been many follow-up research progresses. For example, Sun and Liu \cite{SunWang} discussed the degree theory of the Kazdan-Warner equation on a finite graph. Keller and Schwarz \cite{KellerSchwarz} and Ge and Jiang \cite{GeJiang} studied the same type of equations on a canonically compactifiable graph and an infinite graph, respectively. Existence and multiplicity of solutions for several nonlinear equations on graphs, such as the heat equation, Yamabe equation, Chern-Simons equation etc., were investigated in \cite{Lin1,ZhangLin1,ZhangLin2,Akduman-Pankov,LiuYang,LinYang1,ChaoHou,Wang, HuaXu, HuaLi, ChangWang}. In particular, there were many interesting works for $p$-Laplace equations on graphs. Ge \cite{Ge2,Ge1} studied a $p$-Kazdan-Warner equation and a $p$-Yamabe equation on finite graphs. Peng, Zhu and Zhang \cite{PZZ} considered the existence of multiple solutions to a perturbed $p$-Yamabe equation on a finite graph. For more related results, we refer the readers to \cite{ZhangZhao,Han-Shao,Shao,Shao1} and the references therein.

The main concern of this paper is to discuss the $p$-Laplace equation
\begin{equation}\label{pequation}
-\Delta_pu+\rho|u|^{p-2}u=\psi(x,u)
\end{equation}
on a locally finite graph $G=(V,E)$ for $p>1$,  where $\Delta_{p}$ is the variational $p$-Laplacian on graphs as in \cite{Gri2}, $\rho$ is a potential function and $\psi$ is a nonlinear term. The $p$-Laplacian operator arises in many physical phenomena, such as non-Newtonian fluid, flow through porous media, etc. and kinds of nonlinear $p$-Laplace equations \eqref{pequation} on Euclidean space have been extensively studied, see for examples \cite{Ding,Cao,doo1,doo2,Alv,AdiYang,Yang}. On graphs, the discrete $p$-Laplacian was introduced in \cite{Nakamura} and has been well studied ever since, mostly in the context of nonlinear potential and spectral theory, cf. \cite{Iannizzotto,Chang,Mugnolo,HuaWang} for its historical overviews and recent progresses.

For a $p$-Yamabe type equation on a finite graph, Grigor'yan, Lin and Yang \cite{Gri2} studied the existence of nontrivial solutions by focusing on the variational structure of the equation. In fact, for a finite graph $G=(V,E)$, the dimensions of Sobolev spaces $W^{m,p}(V)$ are all finite and whence the spaces are pre-compact. For this reason, variational problems on finite graphs are comparatively easy to deal with. For a graph $G=(V,E)$ with infinite vertices, Sobolev spaces $W^{m,p}(V)$ are infinite dimensional and their embeddings become unusual. For nonlinear $p$-Laplace equations defined on locally finite graphs, Han and Shao \cite{Han-Shao} proved the existence of positive solutions and Shao \cite{Shao1} investigated the multiplicity of solutions for $p\geq 2$. It is worth pointing out that even if it is just to ensure the reflexivity of the space $W^{1,p}(V)$, the graphs in \cite{Han-Shao,Shao1} (see Corollary 5.8 in \cite{Han-Shao}) are required to satisfy the condition
\begin{equation}\label{M}
  M:=\underset{x\in V}\sup\frac{\underset{y\in V}\sum w_{xy}}{\sigma(x)}<+\infty,
\end{equation}
where $w_{xy}$ is the weight of the edge $xy\in E$ and $\sigma(x)$ is the measure at $x\in V$. Quite recently, we \cite{SYZ} improved the theory of Sobolev spaces on locally finite graphs with the help of a new linear space composed of vector-valued functions with variable dimensions. In particularly, we investigated properties of $W^{m,p}(V)$ for $m\in \mathbb{N}$ and $1\leq p\leq +\infty$ in \cite{SYZ}, including completeness, reflexivity, and separability, on an arbitrary connected locally finite graph without any other restrictions, which are fundamental when dealing with equations on graphs under the variational framework.

Besides existence, for equations with a potential function, the asymptotic behaviour of solutions has also been well studied on graphs. For the special case $p=2$, Zhang and Zhao \cite{ZhangZhao} studied the limit of ground state solutions as $\theta\rightarrow +\infty$ to the following equation
\begin{equation*}
  -\Delta u+(\theta a+1)u=|u|^{q-1}u \quad{\rm in}\ \ \  V,
\end{equation*}
where $q\geq2$, $a(x)\geq 0$ and the potential well $\Omega=\{x\in V:a(x)=0\}$ is a non-empty connected and finite subset of $V$.
While for $p>2$, Han and Shao \cite{Han-Shao} proved a similar result for the following $p$-Laplace equation
\begin{equation*}
  -\Delta_p u+(\theta a+1)|u|^{p-1}u=\psi(x,u) \quad{\rm in}\ \ \  V.
\end{equation*}
In addition, Xu and Zhao \cite{XuZhao} considered the existence and convergence of solutions for a nonlinear elliptic systems and recently, results in \cite{XuZhao} were generalized by Shao \cite{Shao} to a $p$-Laplace system with a homogeneous nonlinearity on a locally finite graph for $p\geq 2$.

In this paper, we consider the existence and convergence of strictly positive solutions of \eqref{pequation}. According to a simple reduction argument (see Sec. \ref{mp}), a strictly positive solution of \eqref{pequation} is equivalent to a nontrivial solution of the following equation
\begin{equation}\label{p-equation}
	-\Delta_pu+\rho|u|^{p-2}u=\psi(x,u^+),
\end{equation}
where $u^+=\max\{u,0\}$ is the positive part of $u$. For this reason, in the remaining part of this paper, the discussion will focus on the equation \eqref{p-equation}.

Compared to the mentioned works about the $p$-Laplace equations, our results have several innovations and improvements. Firstly, the existing works related to $p$-Laplace equations on locally finite graphs all deal with the equation for $p\geq 2$, but we prove the existence results for any $p>1$. Secondly, we aim to study the $p$-Laplace equation \eqref{p-equation} on a locally finite graph $G=(V,E)$ without the condition \eqref{M} on degree and measure of the graph, which is needed in all the mentioned works. Furthermore, when investigate the asymptotic behaviour of solutions to the equation, we consider a more general potential function $\rho=\theta a+b$ instead of $\rho=\theta a+1$, where $b|_\Omega\equiv 0$ is allowed. Finally, we also remove the growth condition of the nonlinear term $\psi(x,s)$ as $s\rightarrow +\infty$, which is used in \cite{Han-Shao, Shao, Shao1}.

To achieve these goals, our methods are also different with those in the existing works. With the help of a discrete gradient for $u$ at $x\in V$, we prove that the energy functional satisfies the $(PS)$ condition in a concise way, which is different from the methods used in \cite{Han-Shao,Shao1}. Moreover, when considering the existence of a ground state solution, \cite{ZhangZhao} used a straightforward approach for the nonlinear term $|u|^{q-1}u$, which is not applicable for our general nonlinear term $\psi(x,u)$. In \cite{Han-Shao}, Han and Shao used a profound deformation lemma to overcome this difficulty. In our paper, different from both \cite{ZhangZhao} and \cite{Han-Shao}, we use an analysis and direct technique inspired by Adimurthi (\cite{Adimurthi}, Lemma 3.5) to get the ground state solution of \eqref{p-equation}.

To state our results in details, let us first introduce some notations and assumptions. Throughout this paper, we always assume that $G=(V,E)$ is a connected and locally finite graph, where $V$ denotes the set of vertices and $E$ denotes the set of edges. The weight $w_{xy}\in \mathbb{R}$ of any edge $xy \in E$ is positive and satisfies $w_{xy}=w_{yx}$. The measure $\sigma$ on $V$ is assumed to have a positive lower bound, i.e. there
exists some positive constant $\sigma_0$ such that
\begin{equation}\label{lowerbd}
	\sigma(x)\geq \sigma_0>0, \quad\forall x\in V.
\end{equation}
A bounded domain $\Omega\subset V$ is always a finite set such that for any $x,y\in\Omega$, the distance $dist(x,y)$ is uniformly bounded from above and the boundary $\partial \Omega$ of $\Omega$ is
$$\partial\Omega=\{y\not\in \Omega:\exists x\in \Omega,\,{\rm such\,\,that}\,\, xy\in E\}.$$
Consequently, we denote $\overline{\Omega}=\Omega\cup\partial\Omega$.

Appropriate concepts of integral and function spaces are also necessary to apply variational methods on graphs. The integral of a function $u$ over $V$ is
\begin{displaymath}
	\int_{V}^{}ud\sigma=\sum_{x\in V}^{}  \sigma(x)u(x).
\end{displaymath}
For any $q>0$, $L^q(V)$ is a linear space of functions $u:V\rightarrow \mathbb{R}$ with the norm
$$
\Vert u\Vert_q=\left(\int_{V}|u|^q d\sigma\right)^{1/q}<+\infty.
$$
While for $q=+\infty$, $L^\infty(V)$ is the space with the norm
$$
\Vert u\Vert_\infty=\sup_V |u|<+\infty.
$$
$W^{1,p}(V)$ is a linear space of functions $u:V\rightarrow\mathbb{R}$ with the norm
$$
\|u\|_{W^{1,p}(V)}=\left(\int_V(|\nabla u|^p+|u|^p)d\sigma\right)^{1/p}<+\infty,
$$
where the gradient of $u$ at $x\in V$ is
\begin{equation}\label{gradient}
\nabla u(x)=\left(\sqrt{\frac{w_{xy_1}}{2\sigma(x)}}(u(y_1)-u(x)),\cdots,\sqrt{\f{w_{xy_{\ell_x}}}{2\sigma(x)}}(u(y_{\ell_x})-u(x))\right),
\end{equation}
and $\{y_1,\cdots,y_{\ell_x}\}$ is the set of all neighbors of $x$ for some positive integer $\ell_x$. Consequently, the norm of $|\nabla u|$ at $x\in V$ is
\begin{displaymath}
	\left | \nabla u \right | (x)=\left(\frac{1}{2\sigma(x)}\sum_{y\sim x}^{}w_{xy} \big (u(y)-u(x)\big )^{2} \right )^{1/2}.
\end{displaymath}
Consequently, $W_0^{1,p}(V)$ is the completion of $C_c(V)$ under the norm $\|\cdot\|_{W^{1,p}(V)}$, where $C_c(V)$ is the set of all functions with finite supports in $V$.
The variational $p$-Laplace operator $\Delta_p$ arises from the variation of $\|\nabla u\|_p^p$ with respect to $u$, or equivalently, it satisfies
$$
\int_V(\Delta_pu)vd\sigma=-\int_V|\nabla u|^{p-2}\nabla u\nabla vd\sigma, \ \ \forall v\in C_c(V),
$$
which also reads as a formula of integration by parts. Point-wisely, 
$$\Delta_pu(x)=\f{1}{\sigma(x)}\sum_{y\sim x}\f{|\nabla u|^{p-2}(x)+|\nabla u|^{p-2}(y)}{2}w_{xy}(u(y)-u(x)),$$
where $y\sim x$ means that $y$ is connected to $x$ by an edge $xy\in E$. If $\Omega\subset V$ is a bounded domain, a straightforward calculation \cite{Han-Shao} shows that
$$
\int_{\Omega}(\Delta_pu)vd\sigma=-\int_{\overline{\Omega}}|\nabla u|^{p-2}\nabla u\nabla vd\sigma, \ \ \forall v\in C_c(\Omega),
$$
where $C_c(\Omega)$ is the set of functions equal to zero on $V\setminus \Omega$.

Let $d(x)=dist(x,O)$ be the distance of $x\in V$ to some fixed vertex $O\in V$. Hereafter we always assume that $\rho:V\rightarrow\mathbb{R}$ is a nonnegative function satisfying
\begin{equation}\label{h-hypothesis}
	\rho(x)\rightarrow+\infty\,\,{\rm as}\,\, d(x)\rightarrow+\infty.
\end{equation}
It follows that for any $\epsilon>0$, there exists some $R_\epsilon>0$ depending only on $\epsilon$ such that
\begin{equation}\label{h-geq}
	\rho(x)\geq \frac{1}{\epsilon},\quad\forall x\in V\setminus B_{R_\epsilon}(O),
\end{equation}
and
\begin{equation}\label{positive-point}
	\rho(x^*)>0\quad{\rm for\,\,\,some}\quad x^*\in B_{R_\epsilon}(O).
\end{equation}
The Banach space $\mathcal{H}_p(V)$, which is the completion of $C_c(V)$ under the norm
\begin{equation}\label{Hp-norm}
	\|u\|_{\mathcal{H}_p(V)}=\left(\int_V(|\nabla u|^p+\rho|u|^p)d\sigma\right)^{1/p}<+\infty,
\end{equation}
is used to  handle the functional related to the equation \eqref{p-equation}. From now on, we may omit the domains of the functions in the Sobolev spaces and just use $L^q$, $\mathcal{H}_p$ etc. for brevity.

Finally, for the nonlinear term $\psi(x,s):V\times [0,+\infty)\rightarrow \mathbb{R}$, we need the following assumptions: \\ [1.2ex]
$(P1)$ $\psi(x,0)=0$ for any $x\in V$ and $\psi(x,s)$ is continuous in $s\in[0,+\infty)$. For any $M>0$, there exists a constant $A_M$ such that $\psi(x,s)\leq A_M$ for $(x,s)\in V\times[0,M]$.\\[1.2ex]
$(P2)$ For any $s>0$ and $x\in V$, there holds
$$0<\alpha \Psi(x,s)\leq s\psi(x,s)$$
for some constant $\alpha>p$ and $\Psi(x,s)=\int_0^s\psi(x,t)dt$ is the primitive of $\psi$.\\[1.2ex]
$(P3)$ There holds
$$\limsup_{s\rightarrow 0+}\f{\psi(x,s)}{s^{p-1}}<\lambda_p=\inf_{u\in \mathcal{H}_p,\,u\not\equiv 0}
\f{\|u\|_{\mathcal{H}_p}^p}{\|u\|_p^p}$$
uniformly in $x\in V$.\\[1.2ex]
$(P4)$ $\psi(x,s)/s^{p-1}$ is strictly increasing in $s>0$ for all $x\in V$.
\vskip6pt
Now we can state our first main result.

\begin{theorem}\label{positive solution} Let $G=(V,E)$ be a connected and locally finite graph with a measure $\sigma$ satisfying \eqref{lowerbd}, $\rho$ be a nonnegative function satisfying \eqref{h-hypothesis}, and suppose that $\psi$ satisfies the hypotheses $(P1)$, $(P2)$ and $(P3)$. Then the equation \eqref{p-equation} has a strictly positive solution in $\mathcal{H}_p(V)$ for any $p>1$.
\end{theorem}

Our second theorem confirms the existence of a strictly positive ground state solution of \eqref{p-equation}.

\begin{theorem}\label{Nehari}Let $G=(V,E)$ be a connected and locally finite graph with a measure $\sigma$ satisfying \eqref{lowerbd}, $\rho$ be a nonnegative function satisfying \eqref{h-hypothesis}, and suppose that $\psi$ satisfies the hypotheses $(P1)$, $(P2)$, $(P3)$ and $(P4)$. Then the equation \eqref{p-equation} has a strictly positive ground state solution in $\mathcal{H}_p(V)$ for any $p>1$. .
\end{theorem}

Next, let us consider a potential $\rho(x)$ in \eqref{p-equation} with the form $\theta a(x)+b(x)$, namely
\begin{equation}\label{potential}
  -\Delta_pu+(\theta a+b)|u|^{p-2}u=\psi(x,u^+)\quad{\rm in}\ \ \  V.
\end{equation}
where $\theta>0$ is the parameter of the potential and $a(x)$ and $b(x)$ are two nonnegative functions. We suppose that $\Omega=\{x\in V:a(x)=0\}$ is a nonempty connected finite subset of $V$, $a(x)\geq 1$ for all $x\in V\setminus\Omega$ and $a(x)\rightarrow+\infty$ as $d(x)\rightarrow+\infty$.
To present the asymptotic behaviour of solutions to the equation \eqref{potential} as $\theta\rightarrow +\infty$, we shall introduce the limit problem which is
defined on the potential well $\Omega$:
\begin{equation}\label{po}\left\{\begin{array}{lll}
	-\Delta_pu+b|u|^{p-2}u=\psi(x,u^+) &{\rm in}&\Omega\\[1.2ex]
	u=0&{\rm on}&\partial\Omega.
\end{array}\right.
	\end{equation}
Let $\mathcal{H}_{a,b}(V)$ be a completion of $C_c(V)$ under the norm
\begin{equation}\label{hab}
\|u\|_{\mathcal{H}_{a,b}}=\left(\int_V(|\nabla u|^p+(a+b)|u|^p)d\sigma\right)^{1/p}<+\infty.
\end{equation}
Under the assumptions $(P1)$, $(P2)$, $(P4)$ and \\[1.2ex]
$(P3^\prime)  \limsup_{s\rightarrow 0+}\frac{\psi(x,s)}{s^{p-1}}<\lambda_{a,b}=\inf_{u\in \mathcal{H}_{a,b},\,u\not\equiv 0}\frac{\|u\|_{\mathcal{H}_{a,b}}^p}{\|u\|_p^p},$
\\[1.2ex]
\noindent by almost the same discussions as in the proof of Theorem \ref{Nehari}, we can conclude that the equation \eqref{potential} has a strictly positive ground state solution $u_\theta\in \mathcal{H}_{a,b}(V)$ for each $\theta>0$. The next theorem tells us what happens as the parameter $\theta$ tends towards $+\infty$.

\begin{theorem}\label{posi-limit}
Suppose that $\psi$ satisfies $(P1)$, $(P2)$, $(P4)$ and $(P3^\prime)$ and $\{u_k\}\subset\mathcal{H}_{a,b}(V)$ is a sequence of positive ground state solutions of the equation \eqref{potential} for $\theta_k\geq 1$ with $\theta_k\rightarrow+\infty$ as $k\rightarrow\infty$. Then up to a subsequence, $\{u_k\}$ converges to some $u_0$ in $W^{1,p}(V)$, where $u_0\equiv 0$ on $V\setminus\Omega$ and the restriction of $u_0$ on $\Omega$ (still denoted by $u_0$) is a positive ground state solution of the Dirichlet problem \eqref{po}.
\end{theorem}

%\begin{remark}
%We are mainly concerned with the existence and asymptotic behaviour of solutions to equations \eqref{p-equation}, \eqref{potential} and \eqref{po}. Our methods are different from and our results are more general than previous works. For example, compared to \cite{Han-Shao}, we relax the assumptions for the nonlinearity $\psi(x,s)$ and present the limit problem \eqref{po} in a more general form. Moreover, we no longer need the restriction \eqref{M} on the graph and our results hold for any $p\in (1,+\infty)$.
%\end{remark}

\begin{remark}
We prove existence and asymptotic behaviour of positive solutions of the equation \eqref{p-equation} under our assumptions $(P1)-(P4)$. If we replace the assumptions with the following $(P1^*)-(P4^*)$, after only minor modifications to our proofs, one can also get the same results for nontrivial solutions of the equation \eqref{pequation} as in \cite{Han-Shao}.\\[1.2ex]
	$(P1^*)$ $\psi(x,0)=0$ for any $x\in V$ and $\psi(x,s)$ is continuous in $s\in\mathbb{R}$. For any $M>0$, there exists a constant $A_M$ such that $|\psi(x,s)|\leq A_M$ for $(x,s)\in V\times[-M,M]$.\\[1.2ex]
	$(P2^*)$ There exists some $\alpha>p$ such that for any $s\in\mathbb{R}\setminus\{0\}$, there holds
	$$0<\alpha \Psi(x,s)\leq s\psi(x,s),\quad\forall x\in V.$$
	$(P3^*)$ There holds
	$$\limsup_{s\rightarrow 0}\frac{|\psi(x,s)|}{|s|^{p-1}}<\lambda_p=\inf_{u\in \mathcal{H}_p,\,u\not\equiv 0}
	\frac{\|u\|_{\mathcal{H}_p}^p}{\|u\|_p^p}$$
	uniformly in $x\in V$.\\[1.2ex]
	$(P4^*)$ $\forall x\in V$, the function $s\mapsto \psi(x,s)/|s|^{p-1}$ is strictly increasing in $s\in\mathbb{R}\setminus\{0\}$.
\end{remark}

We organize the remaining part of the paper as follows. In Section \ref{mp}, we prove the existence of a positive solution  of \eqref{p-equation} by the mountain-pass theorem. In Section \ref{gs}, we use the  method of the Nehari manifold to get a positive ground state solution. Finally, in Section \ref{cgs}, we describe the asymptotic behaviour of the positive ground state solutions of \eqref{potential}.

\section{Positive solution of the mountain-pass type}\label{mp}

In this section, we firstly investigate the embedding properties of the function space $\mathcal{H}_p(V)$ defined as in \eqref{Hp-norm}.

\begin{lemma}\label{embedding}
	For any $1<p<+\infty$, the function space $\mathcal{H}_p(V)$ can be embedded compactly into $L^q(V)$ for all $q$ with $p\leq q\leq +\infty$.
\end{lemma}

\proof Take a bounded sequence $\{u_k\}\subseteq\mathcal{H}_p$, namely
\begin{equation*}\label{bdd}
	\|u_k\|_{\mathcal{H}_p}^p=\int_V(|\nabla u_k|^p+\rho|u_k|^p)d\sigma\leq C.
\end{equation*}
It suffices to prove that for any $q\geq p$, there exists a subsequence (still denoted by $\{u_k\}$) and a function $u_\ast$, 
such that $\{u_k\}$ converges to $u_\ast$ in $L^q(V)$.

First we assume that $p\leq q<+\infty$.
Taking $\epsilon=1$ in \eqref{h-geq} and \eqref{positive-point}, we have
\begin{eqnarray}\nonumber
	\|u_k\|_{W^{1,p}(V)}^p&=&\int_V(|\nabla u_k|^p+|u_k|^p)d\sigma\\\nonumber
	&\leq&\int_{V\setminus B_{R_1}}(|\nabla u_k|^p+\rho|u_k|^p)d\sigma+\int_{B_{R_1}}(|\nabla u_k|^p+|u_k|^p)d\sigma\\\nonumber
	&\leq& \|u_k\|_{\mathcal{H}_p}^p+C\int_{B_{R_1}}(|\nabla u_k|^p+\rho|u_k|^p)d\sigma\\\label{1p}
	&\leq& C,
\end{eqnarray}
where the constant $C$ depends only on $\rho$. In \eqref{1p} we can use the equivalence of the two norms $(\int_{B_{R_1}}(|\nabla u|^p+|u|^p)d\sigma)^{1/p}$ and $(\int_{B_{R_1}}(|\nabla u|^p+\rho|u|^p)d\sigma)^{1/p}$, because $B_{R_1}$ is a connected finite set. According to Corollary 1.2 in \cite{SYZ}, we can deduce from \eqref{1p} that up to a subsequence, $\{u_k\}$ converges weakly in $W_0^{1,p}(V)$ and point-wisely for $x\in V$ to some $u_\ast$.
In particular, we have
$$
\int_V|\nabla u_\ast|^pd\sigma=\lim_{k\rightarrow+\infty}\int_V|\nabla u_\ast|^{p-2}\nabla u_\ast\nabla u_kd\sigma\leq\limsup_{k\rightarrow+\infty}\|\nabla u_k\|_p\|\nabla u_\ast\|_p^{p-1}
$$
and
$$
\int_V\rho|u_\ast|^pd\sigma\leq\limsup_{k\rightarrow+\infty}\int_V\rho|u_k|^pd\sigma,
$$
which immediately lead to $u_\ast\in \mathcal{H}_p$. Moreover, for any $\epsilon>0$,
in view of \eqref{lowerbd} and \eqref{h-geq}, there holds
\begin{eqnarray*}
	\int_{V\setminus B_{R_\epsilon}}|u_k-u_\ast|^qd\sigma&\leq& \|u_k-u_\ast\|_{L^\infty(V\setminus B_{R_\epsilon})}^{q-p}\int_{V\setminus B_{R_\epsilon}}|u_k-u_\ast|^pd\sigma\\
	&\leq&\left(\frac{\epsilon}{\sigma_0}\int_V\rho|u_k-u_\ast|^pd\sigma\right)^{\frac{q}{p}-1}\epsilon\int_{V\setminus B_{R_\epsilon}}\rho|u_k-u_\ast|^pd\sigma\\
	&\leq& C\epsilon^{q/p}
\end{eqnarray*}
for some constant $C$ depending only on $\sigma_0$ and the upper bound of $\|u_k\|_{\mathcal{H}_p}$. On the other hand, up to a subsequence, there holds $\int_{B_{R_\epsilon}}|u_k-u_\ast|^pd\sigma=o_k(1)$, where and in the sequel, $o_k(1)\rightarrow 0$ as $k\rightarrow +\infty$. Therefore, up to a subsequence,
$\{u_k\}$ converges to $u_\ast$ in $L^q(V)$.

For the case $q=+\infty$, by using the same notations, we have
\begin{eqnarray*}
	\|u_k-u_\ast\|_\infty&\leq&\|u_k-u_\ast\|_{L^\infty(V\setminus B_{R_\epsilon})}+\|u_k-u_\ast\|_{L^\infty(B_{R_\epsilon})}\\
	&\leq&\left(\frac{\epsilon}{\sigma_0}\int_{V\setminus B_{R_\epsilon}}\rho|u_k-u_\ast|^pd\sigma\right)^{1/p}+o_k(1)\\
	&\leq&C\epsilon^{1/p}+o_k(1),
\end{eqnarray*}
which implies that up to a subsequence, $\{u_k\}$ converges to $u_\ast$ in $L^\infty(V)$. $\hfill\Box$\\

%\begin{remark}\label{reduction}
Before the proof of Theorem \ref{positive solution}, we point out that nontrvial solutions of \eqref{p-equation}, \eqref{potential} and \eqref{po} must be strictly positive.
%\end{remark}
On a connected graph, this reduction argument is similar to those on Euclidean spaces and can be found in \cite{Gri3} and \cite{Han-Shao}. For the convenience of readers, let us take the equation \eqref{p-equation} as an example to illustrate the details of the argument.
For a solution $u\in \mathcal{H}_p$ of the equation \eqref{p-equation}, let $u^+=\max\{u,0\}$ and $u^-=\min\{u,0\}$ be the positive and negative parts of $u$ respectively. Since
\begin{eqnarray*}\nabla u^+(x)\nabla u^-(x)&=&\f{1}{2\sigma(x)}\sum_{y\sim x}w_{xy}(u^+(y)-u^+(x))(u^-(y)-u^-(x))\\
	&=&-\f{1}{2\sigma(x)}\sum_{y\sim x}w_{xy}(u^+(y)u^-(x)+u^+(x)u^-(y))\\
	&\geq&0,
\end{eqnarray*}
we can easily deduce that
$$\nabla u(x)\nabla u^-(x)=\nabla u^+(x)\nabla u^-(x)+|\nabla u^-|^2(x)\geq |\nabla u^-|^2(x)$$
and
$$|\nabla u|^2(x)=|\nabla u^+|^2(x)+|\nabla u^-|^2(x)+2\nabla u^+(x)\nabla u^-(x)\geq |\nabla u^+|^2(x)+|\nabla u^-|^2(x).$$
Multiplying both sides of \eqref{p-equation} by $u^-$, one has by the above two estimates and integration by parts that
\begin{eqnarray}
	0=\int_V\psi(x,u^+)u^-d\sigma&=&\int_V|\nabla u|^{p-2}\nabla u\nabla u^-d\sigma+\int_V\rho|u|^{p-2}uu^-d\sigma\nonumber\\
	&\geq&\int_V|\nabla u^-|^pd\sigma+\int_V\rho|u^-|^{p}d\sigma\nonumber\\
	&=&\|u^-\|_{\mathcal{H}_p}^p.\label{test}
\end{eqnarray}
Hence $u^-\equiv 0$. This implies that any solution of \eqref{p-equation} would be nonnegative.

Furthermore, suppose that $u(x_0)=0$ for some $x_0\in V$. The equation \eqref{p-equation} implies that $\Delta_pu(x_0)=0$, which leads to
$$
\sum_{y\sim x_0}w_{x_0y}(|\nabla u|^{p-2}(y)+|\nabla u|^{p-2}(x_0))u(y)=0.
$$
As a consequence, $u(y)=0$ for all $y\sim x_0$. Since $V$ is connected, one obtains $u\equiv 0$ on $V$ by repeating this process.
To summarize, a solution $u$ of \eqref{p-equation} satisfies either $u(x)>0$ for all $x\in V$, or $u(x)\equiv 0$ for all $x\in V$.
Then we can conclude that if $u$ is a nontrivial solution of \eqref{p-equation}, in fact, $u$ shall be a strictly positive solution of \eqref{p-equation}. 

In the inequality \eqref{test}, $u^-$ is used as a test function which requires $u^-\in \mathcal{H}_p$ if $u\in \mathcal{H}_p$. On a locally finite graph, this requirement is not clearly met and we present a proof here.

\begin{lemma}
	For any $u\in \mathcal{H}_p(V)$, there exist two sequences $\{\phi_k^{\pm}\}\subset C_c(V)$ such that $$\lim_{k\rightarrow+\infty}\|u^{\pm}-\phi_k^{\pm}\|_{\mathcal{H}_p}=0,$$ namely $u^{\pm}\in \mathcal{H}_p(V)$.  
\end{lemma}
\proof
Because the treatments for the positive and negative parts of $u$ are the same, we will only prove that the positive part $u^+$ can be approximated by a sequence $\{\phi_k^+\}\subset C_c(V)$ in $\mathcal{H}_p(V)$.

Since $u\in\mathcal{H}_p$, we can find a sequence $\{\phi_k\}\subset C_c(V)$ such that as $k\rightarrow+\infty$, $\{\phi_k\}$ converges to $u$ strongly in $\mathcal{H}_p$ and point-wisely in $V$. Moreover, since $\left| (|\nabla \phi_k|-|\nabla u|)(x)\right|\leq |\nabla (\phi_k-u)|(x)$ for any $x\in V$,  there holds for any subset $\tilde{V}\subseteq V$,
\begin{equation}\label{phiku}
	\lim_{k\rightarrow+\infty}\int_{\tilde{V}} |\nabla\phi_k|^p d\sigma= \int_{\tilde{V}} |\nabla u|^p d\sigma
\end{equation}

Obviously for any $R>1$, the positive parts of $u$ and $\phi_k$ satisfy
\begin{equation}\label{phik+}
	\lim_{k\rightarrow+\infty}\int_{B_R} |\nabla(\phi_k^+-u^+)|^p d\sigma=0.
\end{equation}
Since $|\nabla u^+(x)|\leq |\nabla u(x)|$ and $|\nabla \phi_k^+(x)|\leq |\nabla \phi_k(x)|$ for any $x\in V$, we have
\begin{eqnarray*}
	\int_V |\nabla(\phi_k^+-u^+)|^p d\sigma&=&\int_{V\setminus B_R} |\nabla(\phi_k^+-u^+)|^p d\sigma
	+\int_{B_R} |\nabla(\phi_k^+-u^+)|^p d\sigma\\
	&\leq&2^p\int_{V\setminus B_R} |\nabla \phi_k^+|^p d\sigma+2^p \int_{V\setminus B_R} |\nabla u^+|^p d\sigma
	+\int_{B_R} |\nabla(\phi_k^+-u^+)|^p d\sigma\\
	&\leq&2^p\int_{V\setminus B_R} |\nabla \phi_k|^p d\sigma+2^p \int_{V\setminus B_R} |\nabla u|^p d\sigma
	+\int_{B_R} |\nabla(\phi_k^+-u^+)|^p d\sigma\\
	&=&2^{p+1} \int_{V\setminus B_R} |\nabla u|^p d\sigma
	+\int_{B_R} |\nabla(\phi_k^+-u^+)|^p d\sigma+o_k(1),
\end{eqnarray*}
where we have used \eqref{phiku} in the last equality. Then \eqref{phik+} and the arbitrariness of $R$ give us
$$
\lim_{k\rightarrow+\infty}\int_V |\nabla(\phi_k^+-u^+)|^p d\sigma=0.
$$
On the other hand, there obviously holds
$$
\lim_{k\rightarrow+\infty}\int_V \rho|\phi_k^+-u^+|^p d\sigma=0.
$$
Thus the lemma is proved.
$\hfill\Box$\\

To begin the proof of Theorem \ref{positive solution}, for any fixed real number $p>1$, we define a functional $J:\mathcal{H}_p\rightarrow\mathbb{R}$ by
\begin{equation*}\label{functional}
	J(u)=\frac{1}{p}\int_V(|\nabla u|^p+\rho|u|^p)d\sigma-\int_V\Psi(x,u^+)d\sigma.
\end{equation*}
For any $\varphi\in\mathcal{H}_p$, we have
\begin{equation}\label{deriv}\langle J^\prime(u),\varphi\rangle=\int_V|\nabla u|^{p-2}\nabla u\nabla \varphi d\sigma+\int_V\rho|u|^{p-2}u\varphi d\sigma-
	\int_V\psi(x,u^+)\varphi d\sigma.\end{equation}
In particular, if $u\in\mathcal{H}_p$ satisfies $J^\prime(u)=0$, $u$ is a weak solution and also a point-wise solution of the $p$-laplace equation \eqref{p-equation}. Since our method is based on the mountain-pass theorem due to Ambrosetti-Rabinowitz \cite{Ambrosetti-Rabinowitz}, we first prove several lemmas to describe the geometric profile of $J$.

\begin{lemma}\label{lemma1}
There exists some $u\in\mathcal{H}_p$ with $u(x)\geq 0$ for all $x\in V$ such that $J(tu)\rightarrow-\infty$ as $t\rightarrow+\infty$.
\end{lemma}
\proof
By $(P1)$ and $(P2)$, there exist two positive constants $c_1$ and $c_2$ satisfying $\Psi(x,s)\geq c_1s^\alpha-c_2$ for all
$(x,s)\in V\times [0,+\infty)$, where $\alpha>p$ is given as in $(P2)$. For any fixed point $x_0\in V$, we set
$$u(x)=\le\{\begin{array}{lll}
1&{\rm when}& x=x_0\\[1.5ex]
0&{\rm when}& x\not=x_0.
\end{array}\right.$$
Noting that $\alpha>p$, $\sigma(x_0)>0$ and $V$ is locally finite, we calculate for $t>0$,
\begin{eqnarray*}
J(tu)&=&\frac{t^p}{p}\left(\sum_{y\sim x_0}\sigma(y)|\nabla u|^p(y)+\sigma(x_0)|\nabla u|^p(x_0)\right)+\frac{t^p}{p}\sigma(x_0)\rho(x_0)-\sigma(x_0)\Psi(x_0,t)\\
&\leq&\frac{t^p}{p}\left(\sum_{y\sim x_0}\sigma(y)|\nabla u|^p(y)+\sigma(x_0)|\nabla u|^p(x_0)\right)+\frac{t^p}{p}\sigma(x_0)\rho(x_0)-c_1t^\alpha\sigma(x_0)+c_2\sigma(x_0)\\
&\rightarrow&-\infty
\end{eqnarray*}
as $t\rightarrow+\infty$. $\hfill\Box$\\

Secondly we have the following lemma.

\begin{lemma}\label{lemma2}
There exist constants $\delta>0$ and $r>0$ such that $J(u)\geq \delta$ for all $u$ with $\|u\|_{\mathcal{H}_p}=r$.
\end{lemma}

\proof Keep in mind that $\psi(x,s)>0$ for all
$(x,s)\in V\times (0,+\infty)$ and $\psi(x,0)=0$. By $(P3)$, there exist positive constants $\tau$ and $\zeta$ such that if $0<s\leq \zeta$,
there holds $$\Psi(x,s)=\int_0^s\psi(x,t)dt\leq \frac{\lambda_p-\tau}{p}s^p.$$
Note that by $(P2)$, $\Psi(x,s)>0$ for all $s>0$. It follows that
if $s\geq \zeta>0$, for any integer $N>p$, we have
$$\Psi(x,s)\leq \frac{1}{\zeta^N}s^{N}\Psi(x,s).$$
For all $(x,s)\in V\times [0,+\infty)$, there holds
 $$\Psi(x,s)\leq\frac{\lambda_p-\tau}{p}s^p+\frac{1}{\zeta^N}s^N\Psi(x,s).$$
By Lemma \ref{embedding}, for any function $u$ with $\|u\|_{\mathcal{H}_p}\leq 1$,
we have $\|u\|_{\infty}\leq C\|u\|_{\mathcal{H}_p}$, $\|u\|_{N}\leq C\|u\|_{\mathcal{H}_p}$ for some constant $C$ and
$$
\int_Vu^N\Psi(x,u^+)d\sigma\leq \left(\max_{(x,s)\in V\times[0,C]}\Psi(x,s)\right)\,\int_V|u|^Nd\sigma\leq C\|u\|_{\mathcal{H}_p}^N,
$$
where $(P1)$ is employed. Hence we have for any $u$ with $\|u\|_{\mathcal{H}_p}\leq 1$,
\begin{eqnarray*}
  J(u)&\geq& \frac{1}{p}\|u\|_{\mathcal{H}_p}^p-\frac{\lambda_p-\tau}{p}\int_V|u|^pd\sigma-C\|u\|_{\mathcal{H}_p}^N\\
  &\geq& \left(\frac{1}{p}-\frac{\lambda_p-\tau}{p\lambda_p}\right)\|u\|_{\mathcal{H}_p}^p-C\|u\|_{\mathcal{H}_p}^N\\
  &=&\left(\frac{\tau}{p\lambda_p}-C\|u\|_{\mathcal{H}_p}^{N-p}\right)\|u\|_{\mathcal{H}_p}^p.
\end{eqnarray*}
Let $r^{N-p}=\min\{1,{\tau}/(2pC\lambda_p)\}$. We have $J(u)\geq \tau r^p/(2p\lambda_p)$ for all $u$ with $\|u\|_{\mathcal{H}_p}=r$. This gives the desired result. $\hfill\Box$\\

Besides the mountain-pass geometry, to verify the following Palais-Smale condition for the functional $J$ is also important for existence of solutions.

\begin{lemma}\label{ps-cond} The functional $J$ satisfies the  $(PS)_c$ condition for any $c\in\mathbb{R}$ and $p>1$. Namely, if $J(u_k)\rightarrow c$ and $J^\prime(u_k)\rightarrow 0$ for a sequence $\{u_k\}\subset \mathcal{H}_p$, there exists some $u\in\mathcal{H}_p$ such that up to a subsequence, $u_k\rightarrow u$ in $\mathcal{H}_p$.
\end{lemma}

\proof Note that $J(u_k)\rightarrow c$ and $J^\prime(u_k)\rightarrow 0$ as $k\rightarrow+\infty$ are equivalent to
\begin{eqnarray}\label{j-0}
  	\frac{1}{p}\|u_k\|_{\mathcal{H}_p}^p-\int_V\Psi(x,u_k^+)d\sigma=c+o_k(1)
\end{eqnarray}
and
\begin{eqnarray}\label{j1-0}
    \left|\int_V(|\nabla u_k|^{p-2}\nabla u_k\nabla\varphi+\rho|u_k|^{p-2}u_k\varphi-\psi(x,u_k^+)\varphi) d\sigma\right|= o_k(1)\|\varphi\|_{\mathcal{H}_p},
  \quad \forall \varphi\in \mathcal{H}_p.
\end{eqnarray}
Take $\varphi=u_k$ in \eqref{j1-0} and we get
\begin{equation}\label{101}
	\|u_k\|_{\mathcal{H}_p}^p=\int_V\psi(x,u_k^+)u_kd\sigma+o_k(1)\|u_k\|_{\mathcal{H}_p}.
\end{equation}
In view of $(P2)$, we have by combining \eqref{j-0} and \eqref{101} that
\begin{eqnarray*}
  \|u_k\|_{\mathcal{H}_p}^p&=&p\int_V\Psi(x,u_k^+)d\sigma+pc+o_k(1)\\
  &\leq&\f{p}{\alpha}\int_V\psi(x,u_k^+)u_k^+d\sigma+pc+o_k(1)\\
  &=&\f{p}{\alpha}\|u_k\|_{\mathcal{H}_p}^p+o_k(1)\|u_k\|_{\mathcal{H}_p}+pc+o_k(1).
\end{eqnarray*}
Since $\alpha>p$, we know that $u_k$ is bounded in ${\mathcal{H}_p}$, and thus there exists a constant $M$ such that $\|u_k\|_\infty\leq M$. Lemma \ref{embedding} implies that up to a subsequence, $u_k\rightharpoonup u$ weakly in $W_0^{1,p}(V)$ and $u_k\rightarrow u$ in $L^q(V)$ for any $p\leq q\leq+\infty$. It follows that there exists some $\epsilon_0>0$ such that
\begin{eqnarray}\nonumber
  \left|\int_V\psi(x,u_k^+)(u_k-u)d\sigma\right|&\leq&\int_{0\leq |u_k|\leq\epsilon_0}\psi(x,u_k^+)|u_k-u|d\sigma+
  \int_{|u_k|>\epsilon_0}\psi(x,u_k^+)|u_k-u|d\sigma\\\nonumber
  &\leq&\lambda_p\int_{0\leq |u_k|\leq\epsilon_0}u_k^{p-1}|u_k-u|d\sigma+\frac{A_M}{\epsilon^{p-1}}\int_{|u_k|>\epsilon_0}u_k^{p-1}|u_k-u|d\sigma\\\nonumber
  &\leq&\left(\lambda_p+\frac{A_M}{\epsilon_0^{p-1}}\right)\int_V|u_k|^{p-1}|u_k-u|d\sigma\\\nonumber
  &\leq&\left(\lambda_p+\frac{A_M}{\epsilon_0^{p-1}}\right)\|u_k\|_p^{p-1}\|u_k-u\|_p\\\label{f-0}
  &=&o_k(1).
\end{eqnarray}
By subsituting $\varphi=u_k-u\in \mathcal{H}_p$ into \eqref{j1-0}, we obtain by \eqref{f-0} that
\begin{equation}\label{gr-0}
	\int_V\left(|\nabla u_k|^{p-2}\nabla u_k\nabla(u_k-u)+\rho|u_k|^{p-2}u_k(u_k-u)\right)d\sigma=o_k(1).
\end{equation}
Since the weak convergence of $u_k$ in $W_0^{1,p}(V)$ gives
$$
\int_V|\nabla u|^{p-2}\nabla u\nabla (u_k-u)d\sigma=o_k(1)
$$
and H\"older's inequality leads to
\begin{eqnarray*}
	\left|\int_V\rho|u|^{p-2}u(u_k-u)d\sigma\right|&\leq&\left|\int_{V\setminus B_R}\rho|u|^{p-2}u(u_k-u)d\sigma\right|+
	\left|\int_{B_R}\rho|u|^{p-2}u(u_k-u)d\sigma\right|\\&\leq& \left(\int_{V\setminus B_R}\rho|u|^pd\sigma\right)^{1-1/p}\|u_k-u\|_{\mathcal{H}_p}+o_k(1)\\
    &\leq&C\left(\int_{V\setminus B_R}\rho|u|^pd\sigma\right)^{1-1/p}+o_k(1)
\end{eqnarray*}
for any $R>1$, we arrive at
\begin{equation}\label{u-0}
	\int_V\left(|\nabla u|^{p-2}\nabla u\nabla(u_k-u)+\rho|u|^{p-2}u(u_k-u)\right)d\sigma=o_k(1).
\end{equation}

To proceed, we need an elementary inequality which can be found in \cite{Peral} (Lemma A.0.5). 
\begin{equation}\label{inequ}
	\langle|\mathbf{b}|^{p-2}\mathbf{b}-|\mathbf{a}|^{p-2}\mathbf{a},\mathbf{b}-\mathbf{a}\rangle\geq
	 \left\{\begin{array}{lll}
	 	c_p|\mathbf{b}-\mathbf{a}|^p&{\rm if}& p\geq 2\\[1.2ex]
	 	c_p{|\mathbf{b}-\mathbf{a}|^2}{(|\mathbf{b}|+|\mathbf{a}|)^{p-2}}&{\rm if}&1<p<2.
 	\end{array}\right.
 \end{equation}
where $\mathbf{a}$ and $\mathbf{b}$ are two vectors in $\mathbb{R}^n$ with the standard inner product $\langle\cdot,\cdot\rangle$, $c_p$ is a positive constant  depending only on $p$.\\[1.2ex]

{\it Case 1. $p\geq 2$.}

 For any fixed $x\in V$, we know from \eqref{gradient} that $\nabla u(x)$, $\nabla u_k(x)$ and $\nabla (u_k-u)(x)$ are all vectors in $\mathbb{R}^{\ell_x}$. If we take $n=\ell_x$, $\mathbf{b}=\nabla u_k(x)$ and $\mathbf{a}=\nabla u(x)$ in \eqref{inequ}, we immediately obtain that
$$
\int_V|\nabla (u_k-u)|^pd\sigma\leq c_p^{-1}\int_V\left(|\nabla u_k|^{p-2}\nabla u_k-|\nabla u|^{p-2}\nabla u\right)\nabla(u_k-u))d\sigma.
$$
If we take $n=1$, $\mathbf{b}=u_k(x)$ and $\mathbf{a}=u(x)$ in \eqref{inequ}, we obtain that
$$
\int_V\rho|u_k-u|^pd\sigma\leq c_p^{-1}\int_V\rho(|u_k|^{p-2}u_k-|u|^{p-2}u)(u_k-u)d\sigma.
$$
These two estimates together with \eqref{gr-0} and \eqref{u-0} yield
\begin{eqnarray*}
  \|u_k-u\|_{\mathcal{H}_p}^p&\leq&c_p^{-1} \left\{\int_V\left(|\nabla u_k|^{p-2}\nabla u_k\nabla(u_k-u)+\rho|u_k|^{p-2}u_k(u_k-u)\right)d\sigma\right.\\
  &&\quad\quad\left.-\int_V\left(|\nabla u|^{p-2}\nabla u\nabla(u_k-u)+\rho|u|^{p-2}u(u_k-u)\right)d\sigma\right\}\\
  &=&o_k(1).
\end{eqnarray*}
Therefore $u_k\rightarrow u$ in $\mathcal{H}_p$.\\

{\it Case 2. $1<p<2$.}

Since we obviously have $u_k\rightarrow u$ in $\mathcal{H}_p$ for $u_k\equiv u\equiv0$, we will assume that $u_{k}\not\equiv 0$ or $u\not\equiv 0$ next.
For $1<p<2$, it follows from H\"older's inequality that
\begin{equation}\label{1<p<2}
  \int_V|\nabla u_k-\nabla u|^pd\sigma\leq\left(\int_V\frac{|\nabla u_k-\nabla u|^2}{(|\nabla u_k|+|\nabla u|)^{2-p}}d\sigma\right)^{\frac{p}{2}}
  \left(\int_V(|\nabla u_k|+|\nabla u|)^{p}d\sigma\right)^{1-\frac{p}{2}}
\end{equation}
and
\begin{equation}\label{hu}
  \int_V\rho|u_k-u|^pd\sigma\leq\left(\int_V\frac{\rho|u_k-u|^2}{(|u_k|+|u|)^{2-p}}d\sigma\right)^{\frac{p}{2}}
  \left(\int_V\rho(|u_k|+|u|)^{p}d\sigma\right)^{1-\frac{p}{2}}.
\end{equation}
Since $\{u_k\}$ is bounded in $\mathcal{H}_p$, in view of \eqref{inequ}, \eqref{1<p<2} and \eqref{hu}, we have
\begin{eqnarray*}
  \|u_k-u\|_{\mathcal{H}_p}^2&\leq&C\int_V\left(\frac{|\nabla u_k-\nabla u|^2}{(|\nabla u_k|+|\nabla u|)^{2-p}}+
  \frac{\rho|u_k-u|^2}{(|u_k|+|u|)^{2-p}}\right)d\sigma\\
  &\leq&C\int_V\left((|\nabla u_k|^{p-2}\nabla u_k-|\nabla u|^{p-2}\nabla u)\nabla (u_k-u)\right.\\
  &&\quad\left.
  +\rho(|u_k|^{p-2}u_k-|u|^{p-2}u)(u_k-u)\right)d\sigma,
\end{eqnarray*}
which together with \eqref{gr-0} and \eqref{u-0} implies that $\|u_k-u\|_{\mathcal{H}_p}=o_k(1)$. This ends the proof of the lemma. $\hfill\Box$\\

{\it Completion of the proof of Theorem \ref{positive solution}.} Lemmas \ref{lemma1}, \ref{lemma2} and \ref{ps-cond} ensure all the assumptions of the mountain-pass theorem and the Palais-Smale condition for $J$. Consequently we can conclude that  $$c=\min_{\gamma\in\Gamma}\max_{u\in\gamma}J(u)$$
is a critical level of $J$, where
$$
\Gamma=\left\{\g\in C([0,1],\mathcal{H}_p): \g(0)=0, \g(1)=u^\ast\right\}.
$$
In particular, there exists a weak solution $u\in \mathcal{H}_p$ of \eqref{p-equation} such that $J(u)=c$. Moreover, since
$$
J(u)=c\geq \delta>0,
$$
we know that $u\not\equiv 0$ and the proof of the theorem is completed. $\hfill\Box$

\section{Positive ground state solution}\label{gs}

If we can obtain a nontrivial ground state solution of the equation \eqref{p-equation}, the solution must be strictly positive and Theorem \ref{Nehari} can be proved. In this section, we will use the method of Nehari manifold to achieve this aim.

The Nehari manifold associated to the functional $J$ is defined as follows
$$\mathcal{N}_{\psi}=\left\{u\in\mathcal{H}_p\setminus\{0\}: \langle J^\prime(u),u\rangle=0\right\}.$$
We use
$$m=\inf_{u\in\mathcal{N}_{\psi}}J(u).$$
to denote the least energy of the functional among the Nehari manifold. In view of \eqref{deriv}, one has $\|u\|_{\mathcal{H}_p}^p=\int_V\psi(x,u^+)u^+d\sigma$ for any $u\in \mathcal{N}_{\psi}$ and immediately we have  $u^+\not\equiv 0$.\\

We first present a monotonicity lemma which will be frequently used in the following.

\begin{lemma}\label{monotone}
For any function $u$ with $u^+\not\equiv 0$, the function
\begin{equation}\label{g-2}\gamma_u(t)=\frac{1}{t^{p-1}}\int_V\psi(x,tu^+(x))u^+(x)d\sigma\end{equation}
is strictly increasing in $t\in(0,+\infty)$.
\end{lemma}

\proof  For any $x\in V$, if $u(x)>0$, by the hypothesis $(P4)$, the function $\mu(t)=\frac{1}{t^{p-1}}\psi(x,tu^+(x))u^+(x)$ is strictly increasing in $t\in(0,+\infty)$. Otherwise, $u^+(x)=0$, and whence $\mu(t)\equiv 0$ for all  $t\in(0,+\infty)$. Combining these two cases, one completes the proof of this lemma.$\hfill\Box$\\

The next lemma tells us that $\forall u\in\mathcal{H}_p\setminus\{0\}$ with $u^+\not\equiv 0$, a line passing through $u$ and $v\equiv 0$ must pass through the Nehari manifold.

\begin{lemma}\label{unique}
For any $u\in\mathcal{H}_p\setminus\{0\}$ with $u^+\not\equiv 0$, there exists a unique $t_0\in(0,+\infty)$ such that $t_0u\in\mathcal{N}_{\psi}$
and $J(t_0u)=\max_{t\in(0,+\infty)}J(tu)$.
\end{lemma}

\proof For any $t>0$, there holds
\begin{equation*}\label{derivative-2}
	\frac{d}{dt} J(tu)=\langle J^\prime(tu),u\rangle=t^{p-1}\left(\|u\|_{\mathcal{H}_p}^p-
\gamma_u(t)\right),
\end{equation*}
where $\gamma_u(t)$ is defined as in \eqref{g-2}.
In view of $(P3)$, there exists some $\delta>0$ and $0<\tau<\lambda_p$ such that for all $(x,s)\in V\times(0,\delta)$,
\begin{equation}\label{F-est}0\leq \Psi(x,s)\leq \f{\tau}{p}s^p.\end{equation}
Since by Lemma \ref{embedding}, there exists a constant ${c}_0>0$ such that $\|u\|_\infty\leq {c}_0\|u\|_{\mathcal{H}_p}$, we can find some $t_0>0$ such that $tu^+(x)<\delta$ for $0<t<t_0$ and $x\in V$. This together with \eqref{F-est} leads to
\begin{eqnarray}\nonumber
J(tu)&=&\f{1}{p}\|tu\|_{\mathcal{H}_p}^p-\int_V\Psi(x,tu^+)d\sigma\\\nonumber
&\geq&\f{t^p}{p}\|u\|_{\mathcal{H}_p}^p-\f{t^p}{p}\tau\|u^+\|_p^p\\\label{J-tu-2}
&\geq&\f{t^p}{p}\left(1-\f{\tau}{\lambda_p}\right)\|u\|_{\mathcal{H}_p}^p>0
\end{eqnarray}
for $t\in(0,t_0)$. Combining $(P1)$ and $(P2)$, one easily finds $\Psi(x,s)\geq c(s^\alpha-s^p)$ for some positive constant $c$ and
 all $(x,s)\in V\times[0,+\infty)$. Since $\alpha>p$ and $u^+\not\equiv 0$, we obtain
\begin{eqnarray}\nonumber
J(tu)&\leq&\f{t^p}{p}\|u\|_{\mathcal{H}_p}^p-ct^\alpha\int_V(u^+)^\alpha d\sigma+ct^p\int_V(u^+)^pd\sigma\\
&\rightarrow&-\infty\label{-wuq-2}
\end{eqnarray}
as $t\rightarrow+\infty$. It follows from \eqref{J-tu-2}, \eqref{-wuq-2} and $J(0)=0$ that $\mathcal{J}_u(t)=J(tu)$ achieves its maximum
at some $t_0\in(0,+\infty)$. Consequently we have $\mathcal{J}_u^\prime(t_0)=0$ and $t_0u\in\mathcal{N}_{\psi}$. In view of \eqref{derivative-2} and Lemma \ref{monotone}, we conclude that $t_0$ is unique.$\hfill\Box$\\

To search for a ground state solution, we shall give several lemmas about the least energy of $J$ in $\mathcal{N}_{\psi}$.

\begin{lemma}\label{inf}
$m=\inf_{u\in\mathcal{N}_{\psi}}J(u)>0$.
\end{lemma}

\proof $\forall u\in\mathcal{N}_{\psi}$, there holds
\begin{equation}\label{neh-2}
	\|u\|_{\mathcal{H}_p}^p=\int_V\psi(x,u^+)u^+d\sigma.
\end{equation}
In view of $(P3)$, there holds
$$0\leq \psi(x,s)\leq \tau s^{p-1},\quad\forall (x,s)\in V\times [0,\delta)$$
for some constants $\delta>0$ and $0<\tau<\lambda_p$. By Lemma \ref{embedding}, there exists a constant ${c}_0>0$ such that
$\|u\|_\infty\leq {c}_0\|u\|_{\mathcal{H}_p}$. If $u$ satisfies $\|u\|_{\mathcal{H}_p}<\delta/c_0$, we have
$$\int_V\psi(x,u^+)u^+d\sigma\leq \tau\int_V(u^+)^pd\sigma\leq \f{\tau}{\lambda_p}\|u\|_{\mathcal{H}_p}^p,$$
which together with \eqref{neh-2} leads to
$\|u\|_{\mathcal{H}_p}^p\leq ({\tau}/{\lambda_p})\|u\|_{\mathcal{H}_p}^p$.
This is impossible since $u\not\equiv 0$. Hence for any $u\in\mathcal{N}_{\psi}$, there must hold
\begin{equation}\label{uhpbound}
	\|u\|_{\mathcal{H}_p}\geq \delta/c_0.
\end{equation}
Thus by $(P2)$,
\begin{eqnarray*}\nonumber
J(u)
\geq\left(\f{1}{p}-\frac{1}{\alpha}\right)\|u\|_{\mathcal{H}_p}^p
\geq\left(\f{1}{p}-\frac{1}{\alpha}\right)\left(\frac{\delta}{c_0}\right)^p.
\end{eqnarray*}
This implies $m>0$. $\hfill\Box$\\

\begin{lemma}\label{achieved}
 $m=\inf_{u\in\mathcal{N}_{\psi}}J(u)$ is achieved by some $u_*\in\mathcal{N}_{\psi}$.
\end{lemma}
\proof Suppose that $\{u_k\}$ is a sequence in $\mathcal{N}_{\psi}$ such that $J(u_k)\rightarrow m$ as $k\rightarrow\infty$.
Since
$$J(u_k)\geq\left(\f{1}{p}-\f{1}{\alpha}\right)\|u_k\|_{\mathcal{H}_p}^p$$
and $\alpha>p$, $\{u_k\}$ must be a bounded sequence in $\mathcal{H}_p$. Lemma \ref{embedding} confirms that there exists some $u_*\in \mathcal{H}_p$ such that
up to a subsequence, $\{u_k\}$ converges to $u_*$ weakly in $W_0^{1,p}(V)$ and strongly in $L^q(V)$ for all $p\leq q\leq +\infty$. In particular, we have $\|u_k\|_{\infty}$ is uniformly bounded by some positive constant $M$. This together with $(P1)$, $(P3)$ and \eqref{uhpbound} gives
\begin{eqnarray*}
	\left(\frac{\delta}{c_0}\right)^p\leq\|u_k\|_{\mathcal{H}_p}^p&=&\int_{u_k^+< \delta}\psi(x,u_k^+)u_k^+ d\sigma+\int_{M\geq u_k^+\geq \delta}\psi(x,u_k^+)u_k^+ d\sigma\\
	&\leq&\tau\int_{u_k^+< \delta} |u_k^+|^p d\sigma+A_M\int_{M\geq u_k^+\geq \delta}|u_k^+| d\sigma\\
	&\leq&\tau\int_{u_k^+< \delta} |u_k^+|^p d\sigma+\frac{A_M}{\delta^{p-1}}\int_{M\geq u_k^+\geq \delta}|u_k^+|^p d\sigma\\
	&\leq&C\int_V |u_k^+|^p d\sigma=C\int_V |u_*^+|^p d\sigma +o_k(1).
\end{eqnarray*}
where $A_M$ is a constant determined by $(P1)$, and $\tau<\lambda_p$ and $\delta>0$ are constants determined by $(P3)$.
The above fact implies that $u_*^+\not\equiv 0$. As an analog of \eqref{f-0}, we have
\begin{equation}\label{leq-2}
\|u_*\|_{\mathcal{H}_p}^p\leq\lim_{k\rightarrow\infty}\|u_k\|_{\mathcal{H}_p}^p=\lim_{k\rightarrow\infty}\int_V\psi(x,u^+_k)u^+_kd\sigma=
\int_V\psi(x,u_*^+)u_*^+d\sigma.
\end{equation} Moreover, by $(P3)$, there exist two constants $c_1>0$ and $\delta>0$
such that $|\Psi(x,t_1)-\Psi(x,t_2))|\leq c_1(|t_1|^{p-1}+|t_2|^{p-1})|t_1-t_2|$ for all $t_1,t_2\in [0,\delta)$ and all $x\in V$. This
together with H\"older's inequality and the inequality $|a^+-b^+|\leq |a-b|$, $\forall a,b\in\mathbb{R}$, gives us
\begin{equation}\label{small-2}
	\left|\int_{u_k^+<\delta,\,u_*^+<\delta}\left(\Psi(x,u_k^+)-\Psi(x,u_*^+)\right)d\sigma\right|\leq C\left(\int_V|u_k-u_*|^pd\sigma\right)^{1/p}.
\end{equation}
Since $|u_k(x)|+|u_*(x)|\leq c_2$ for all $x\in V$, $k\in\mathbb{N}^\ast$ and some constant $c_2$, and 
$|\Psi(x,t_1)-\Psi(x,t_2)|\leq A_{c_2}|t_1-t_2|$ for all $t_1,t_2\in[0,c_2]$, where $A_{c_2}$ is determined by $(P1)$, we have
\begin{eqnarray*}
&&\left|\int_{V\setminus\left\{x\in V:u_k^+(x)<\delta,\,u_*^+(x)<\delta\right\}}\left(\Psi(x,u_k^+)-\Psi(x,u_*^+)\right)d\sigma\right|\\&&\leq
C\int_{u_k^+\geq \delta}|u_k-u_*|d\sigma+C\int_{u_*^+\geq \delta}|u_k-u_*|d\sigma\\
&&\leq\frac{C}{{\delta}^{p-1}}\int_{u_k^+\geq \delta}|u_k|^{p-1}|u_k-u_*|d\sigma+\frac{C}{{\delta}^{p-1}}\int_{u_*^+\geq \delta}|u_*|^{p-1}|u_k-u_*|d\sigma\\
&&\leq\frac{C}{{\delta}^{p-1}}(\|u_k\|_p^{p-1}+\|u_*\|_p^{p-1})\|u_k-u_*\|_p,
\end{eqnarray*}
where we have used the fact that
$$
V\setminus\left\{x\in V:u_k^+(x)<\delta,\,u_*^+(x)<\delta\right\}=\{x\in V:u_k^+(x)\geq \delta\}\cup\{x\in V:u_*^+(x)\geq\delta\}.
$$
This together with \eqref{small-2} leads to
\begin{equation}\label{converge-2}
	\lim_{k\rightarrow\infty}\int_V\Psi(x,u_k^+)d\sigma=\int_V\Psi(x,u_*^+)d\sigma.
\end{equation}
Combining \eqref{leq-2} and \eqref{converge-2}, we obtain
\begin{equation}\label{energy-2}
J(u_*)\leq \lim_{k\rightarrow\infty}J(u_k)=m.
\end{equation}
There is only left to prove $u_*\in\mathcal{N}_{\psi}$. Suppose not, in view of \eqref{leq-2}, we have
\begin{equation}\label{<-2}
	\|u_*\|_{\mathcal{H}_p}^p<\int_V\psi(x,u_*^+)u_*^+d\sigma.
\end{equation}
Since $u_*^+\not\equiv 0$, by Lemma \ref{unique}, there exists a unique $t_0\in(0,+\infty)$ such that $t_0u_*\in\mathcal{N}_{\psi}$ and $J(t_0u_*)=\max_{t\in(0,+\infty)}J(tu)$. At the same time, since $u_k\in\mathcal{N}_{\psi}$, we also have
$J(u_k)=\max_{t\in(0,+\infty)}J(tu_k)$.  Then it follows from \eqref{leq-2}, \eqref{converge-2} and \eqref{<-2} that
\begin{eqnarray*}
m\leq J(t_0u_*)&=&\frac{t_0^p}{p}\|u_*\|_{\mathcal{H}_p}^p-\int_V\Psi(x,t_0u_*^+)d\sigma\\
&<&\frac{t_0^p}{p}\int_V\psi(x,u_*^+)u_*^+d\sigma-\int_V\Psi(x,t_0u_*^+)d\sigma\\
&=&\lim_{k\rightarrow+\infty}\left(\frac{t_0^p}{p}\int_V\psi(x,u_k^+)u_k^+d\sigma-\int_V\Psi(x,t_0u_k^+)d\sigma\right)\\
&=&\lim_{k\rightarrow+\infty}\left(\frac{1}{p}\|t_0u_k\|_{\mathcal{H}_p}^p-\int_V\Psi(x,t_0u_k^+)d\sigma\right)\\
&\leq&\lim_{k\rightarrow+\infty}J(u_k)=m,
\end{eqnarray*}
which is a contradiction. Therefore $\|u_*\|_{\mathcal{H}_p}^p=\int_V\psi(x,u_*^+)u_*^+d\sigma$, or equivalently $u_*\in\mathcal{N}_{\psi}$. This together with \eqref{energy-2}
gives $J(u_*)=m$.$\hfill\Box$\\

\begin{lemma}\label{critical}
	If $v\in\mathcal{N}_{\psi}$ satisfies $J^\prime(v)\not=0$, there holds $J(v)>m=\inf_{u\in\mathcal{N}_{\psi}}J(u)$.
\end{lemma}

\proof Here we use a modified argument from Adimurthi (\cite{Adimurthi}, Lemma 3.5) to prove the lemma.

Since $J^\prime(v)\not=0$, we can take some $v_0\in \mathcal{H}_p$ such that
$\langle J^\prime(v),v_0\rangle=-1$. For $s,t\in\mathbb{R}$, define $\eta(s,t)=sv+tv_0$. Clearly we have
$$\lim_{(s,t)\rightarrow(1,0)}\frac{d}{dt}J(\eta(s,t))=\langle J^\prime(v),v_0\rangle=-1.$$
Hence there exist two sufficiently small constants $\delta_1>0$ and $\delta_2>0$ such that for all $(s,t)\in[1-\delta_1,1+\delta_1]
\times (0,\delta_2]$, there holds
\begin{equation}\label{str-2}J(\eta(s,t))<J(\eta(s,0))=J(sv).\end{equation}
Define
$$\Upsilon(s,t)=\|\eta(s,t)\|_{\mathcal{H}_p}^p-\int_V\psi(x,(\eta(s,t))^+)(\eta(s,t))^+d\sigma.$$
For any $s\in[1-\delta_1,1+\delta_1]$, there holds
\begin{eqnarray*}
\Upsilon(s,0)&=&s^p\|v\|_{\mathcal{H}_p}^p-\int_V\psi(x,sv^+)sv^+d\sigma\\
&=&s^p\left(\|v\|_{\mathcal{H}_p}^p-\gamma_{v}(s)\right),
\end{eqnarray*}
where $\gamma_{v}(s)$ is defined as in \eqref{g-2}. Since $v^+\not\equiv 0$,  by Lemma \ref{monotone}, $\gamma_{v}(s)$ is strictly increasing in
$s\in(0,+\infty)$. This together with the fact $v\in\mathcal{N}_{\psi}$ implies $\Upsilon(1+\delta_1,0)<\Upsilon(1,0)=0<\Upsilon(1-\delta_1,0)$.
Since $\Upsilon(s,t)$ is continuous in $t\in[0,+\infty)$, there exists a constant $\delta_3\in (0, \delta_2)$ such that
$\Upsilon(1+\delta_1,t)<0<\Upsilon(1-\delta_1,t)$ for all $t\in(0,\delta_3]$. Again by continuity of of $\Upsilon(s,t)$, there exists some $s_t\in (1-\delta_1,1+\delta_1)$ for any $t\in(0,\delta_3]$ such that $\Upsilon(s_t,t)=0$. Since both $v\not\equiv 0$ and $v_{0}\not\equiv0$, it is easy to conclude for $t$ small enough that $\eta(s_t,t)=s_tv+tv_0\not\equiv0$ and $\eta(s_t,t)=s_tv+tv_0\in\mathcal{N}_{\psi}$. Finally, by Lemma \ref{unique} and \eqref{str-2}, we obtain
$$m\leq J(\eta(s_t,t))<J(\eta(s_t,0))=J(s_tv)\leq J(v),$$
and the lemma is proved.$\hfill\Box$\\

{\it Completion of the proof of Theorem \ref{Nehari}.}  Take $u_*\in\mathcal{N}_{\psi}$ given by Lemma \ref{achieved} such that
$J(u_*)=\inf_{u\in\mathcal{N}_{\psi}}J(u)$. Lemma \ref{critical} tells us that $J^\prime(u_*)=0$. Thus $u_*$ is a positive
ground state solution of the $p$-laplace equation \eqref{p-equation}.  $\hfill\Box$

\section{Convergence of positive ground state solutions}\label{cgs}

In this section, we turn to the asymptotic behaviour of solutions of \eqref{potential}.
To be specific, we shall prove that a sequence of positive ground state solutions $\{u_k\}$ of \eqref{potential} with the parameter $\theta_k\rightarrow +\infty$ as $k\rightarrow+\infty$, will converge up to a subsequence to a positive ground state solution of the limit problem \eqref{po}.

To adapt to the special form of the potential $\rho=\theta_k a+b$, we first introduce some notations.
Let $\mathcal{H}_k(V)$ be a completion of $C_c(V)$ under the norm
$$\|u\|_{\mathcal{H}_k}=\left(\int_V(|\nabla u|^p+(\theta_ka+b)|u|^p)d\sigma\right)^{1/p}<+\infty.$$
The functional $J_k: \mathcal{H}_k\rightarrow \mathbb{R}$ of \eqref{potential} is defined as
$$J_k(u)=\f{1}{p}\|u\|_{\mathcal{H}_k}^p-\int_V\Psi(x,u^+)d\sigma.$$
The corresponding Nehari manifold is
$$\mathcal{N}_{k}=\left\{u\in\mathcal{H}_k\setminus\{0\}:\|u\|_{\mathcal{H}_k}^p=\int_V\psi(x,u^+)u^+d\sigma\right\}$$
and the least energy of $J_k$ in the manifold is
$$m_{k}=\inf_{u\in\mathcal{N}_{k}}J_k(u).$$
The suitable working space of \eqref{po} is the Sobolev space $W^{1,p}_{0}(\Omega)$ with the norm
$$\|u\|_{\Omega}=\left(\int_{\overline{\Omega}}|\nabla u|^pd\sigma+\int_{\Omega}b|u|^p)d\sigma\right)^{1/p}<+\infty,$$
where $\overline{\Omega}=\Omega\cup\partial\Omega$.
It is worth pointing out that, since $\Omega$ is connected and $u|_{\partial \Omega}\equiv0$ for any $u\in W_0^{1,p}(\Omega)$, $\|\cdot\|_{\Omega}$ is still a norm even if $b|_{\Omega}\equiv 0$.

We define the functional $J_\Omega: W^{1,p}_{0}(\Omega)\rightarrow \mathbb{R}$ related to \eqref{po} as
$$J_\Omega(u)=\f{1}{p}\|u\|_\Omega^p-\int_\Omega \Psi(x,u^+)d\sigma.$$
The Nehari manifold related to \eqref{po} is defined as
$$\mathcal{N}_{\Omega}=\left\{u\in W_0^{1,p}(\Omega)\setminus\{0\}:\|u\|_\Omega=\int_\Omega \psi(x,u^+)u^+d\sigma\right\}$$
and the least energy of $J_{\Omega}$ in the manifold is
$$m_{\Omega}=\inf_{u\in\mathcal{N}_{\Omega}}J_\Omega(u).$$\\

{\it Completion of the proof of Theorem \ref{posi-limit}.} Since $J_k(u)=J_\Omega(u)$ for all $u\in W_0^{1,p}(\Omega)$ and $\mathcal{N}_{\Omega}\subset\mathcal{N}_{k}$, we have
\begin{equation*}
	m_{k}=\inf_{u\in\mathcal{N}_{k}}J_k(u)\leq m_{\Omega}=\inf_{u\in\mathcal{H}_\Omega}J_\Omega(u)
\end{equation*}
for all $k=1,2,\cdots$. By Theorem \ref{Nehari}, there exists a positive ground state solution $u_k\in\mathcal{N}_{k}$ of \eqref{potential} such that
$J_k(u_k)=m_{k}$ and $J_k^\prime(u_k)=0$. By $(P2)$, $u_k$ satisfies
\begin{eqnarray*}
m_{\Omega}\geq m_{k}&=&J_k(u_k)=\f{1}{p}\|u_k\|_{\mathcal{H}_{k}}^p-\int_V\Psi(x,u_k^+)d\sigma\\
&\geq&\f{1}{p}\|u_k\|_{\mathcal{H}_k}^p-\f{1}{\alpha}\int_V\psi(x,u_k^+)u_k^+d\sigma\\
&=&\left(\f{1}{p}-\f{1}{\alpha}\right)\|u_k\|_{\mathcal{H}_k}^p.
\end{eqnarray*}
Hence one has
$\|u_k\|_{\mathcal{H}_{a,b}}\leq \|u_k\|_{\mathcal{H}_k}\leq C$, where $\|\cdot\|_{\mathcal{H}_{a,b}}$ is defined in \eqref{hab}. Consequently, by Lemma \ref{embedding}, we know that up to a subsequence, $\{u_k\}$ converges to some $u_0$ strongly in $L^q(V)$ for any $p\leq q\leq +\infty$ and point-wisely in $V$. Moreover, we also have $\{u_k\}$ is bounded in $W_0^{1,p}(V)$. Since $W_0^{1,p}(V)$ is reflexive (Corollary 1.2 in \cite{SYZ}), there exists a function
$\tilde{u}_0\in W_0^{1,p}(V)$ such that up to a subsequence, $\{u_k\}$ converges to $\tilde{u}_0$ weakly in $W_0^{1,p}(V)$. The point-wisely convergence of $u_k$ to $u_0$ tells us that in fact, we have $u_0$ and $\tilde{u}_0$ are the same function in $V$. Obviously, $u_{0}(x)\geq 0$ for all $x\in V$. We claim that
\begin{equation}\label{u-equiv-0}u_0\equiv 0\quad{\rm on}\quad V\setminus\Omega.\end{equation}
Otherwise, suppose that $u_0(x_0)>0$ for some $x_0\in V\setminus\Omega$. Since $u_k$ is a solution of \eqref{potential}, at $x_0$ it satisfies
\begin{equation}\label{vari-0}
	(\theta_ka(x_0)+b(x_{0}))u_k(x_0)^{p-1}=\Delta_pu_k(x_0)+\psi(x_0,u_k(x_0)).
\end{equation}
As $k\rightarrow+\infty$, the left side of \eqref{vari-0} tends to infinity, but the right side of \eqref{vari-0} remains bounded, which is a contradiction. This confirms the claim \eqref{u-equiv-0}.

By $(P3^\prime)$ and a discussion as in \eqref{f-0}, there exists two constants $\tau<\lambda_{a,b}$ and $C$ such that
$$\int_V\psi(x,u_k^+)u_k^+ d\sigma\leq \tau\int_Vu_k^pd\sigma+C\int_Vu_k^{p+1}d\sigma.$$
This together with Lemma \ref{embedding} gives
\begin{eqnarray*}
0=\langle J_k^\prime(u_k),u_k\rangle&=&\|u_k\|_{\mathcal{H}_{k}}^p-\int_V\psi(x,u_k^+)u_k^+d\sigma\\
&\geq&\left(1-\f{\tau}{\lambda_{a,b}}\right)\|u_k\|_{\mathcal{H}_{a,b}}^p-C\|u_k\|_{\mathcal{H}_{a,b}}^{p+1}.
\end{eqnarray*}
Hence $\|u_k\|_{\mathcal{H}_{a,b}}\geq \epsilon_0$ for some constant $\epsilon_0>0$, which leads to
$u_0\in W_0^{1,p}(\Omega)\setminus\{0\}$. Moreover, point-wisely convergence of $u_k$ to $u_0$ confirms that $u_0$ satisfies the equation \eqref{po}, namely
\begin{equation*}\left\{\begin{array}{lll}
-\Delta_pu_0+b|u_0|^{p-2}u_0=\psi(x,u_0^+)&{\rm in}&\Omega\\[1.2ex]
u_0=0&{\rm on}&\p\Omega.
\end{array}\right.
\end{equation*}
Multiplying both sides of the above equation with $u_0$, we conclude that $u_0\in \mathcal{N}_{\Omega}$ by integration by parts.

We now prove
\begin{equation}\label{-0-0}
	\lim_{k\rightarrow\infty}\int_V\theta_ka|u_k|^pd\sigma=0\quad{\rm and}\quad \lim_{k\rightarrow\infty}\int_V|\nabla u_k|^pd\sigma=\int_{\overline{\Omega}}|\nabla u_0|^pd\sigma.
\end{equation}
Otherwise, at least one of the two equalities in \eqref{-0-0} fails. As a consequence,
\begin{eqnarray*}
\int_{\overline{\Omega}}|\nabla u_0|^pd\sigma+\int_{\Omega}b|u_0|^pd\sigma
&<&\lim_{k\rightarrow\infty}\int_V(|\nabla u_k|^p+(\theta_ka+b)|u_k|^p)d\sigma\\
&=&\lim_{k\rightarrow\infty}\int_V\psi(x,u_k^+)u_k^+d\sigma\\
&=&\int_V\psi(x,u_0^+)u_0^+ d\sigma,
\end{eqnarray*}
which contradicts $u_0\in\mathcal{N}_\Omega$. Here in the last equality, an analog of \eqref{leq-2} has been used.
Hence \eqref{-0-0} holds.

The power of \eqref{-0-0} is evident. On one hand, it together with Lemma \ref{embedding} gives $\|u_k-u_0\|_{W^{1,p}(V)}=o_k(1)$.
On the other hand, it together with $\mathcal{N}_\Omega\subset\mathcal{N}_k$ leads to
$$J_\Omega(u_0)=\lim_{k\rightarrow\infty}J_k(u_k)=\lim_{k\rightarrow\infty}\inf_{u\in\mathcal{N}_{k}}J_k(u)\leq \inf_{u\in\mathcal{N}_{\Omega}}J_\Omega(u)
=m_{\Omega}.$$
Therefore $u_0$ is a ground state solution of \eqref{po}. Since $u_0(x)\geq 0$ for all $x\in \Omega$ and $J_\Omega(u_0)=m_{\Omega}\geq m_k>0$, we know that $u_0\not\equiv 0$ in $\Omega$. Furthermore, we claim that $u_0(x)>0$ for all $x\in \Omega$. Therefore, $u_0$ is a positive ground state solution of \eqref{po} and the proof of the theorem is completed. $\hfill\Box$\\

\begin{remark}
The proof of the above convergence theorem also implies the existence of a positive ground state solution of the Dirichlet problem \eqref{po}. In Euclidean space, usual variational methods to prove the existence of solutions to nonlinear problems defined on a bounded domain including minimizing the corresponding functional, mountain-pass theorem, etc. Instead of these methods, we directly obtain a ground state solution for the equation \eqref{po} by proving convergence of the ground state solution sequence $\{u_k\}$ of \eqref{potential} in $W^{1,p}(V)$, which is a discussion different from those in Euclidean space. 
\end{remark}

\section*{Acknowledgements}
This research is supported by National Natural Science Foundation of China (No. 12271039 and No. 12101355) and the Open Project Program of Key Laboratory of Mathematics and Complex System, Beijing Normal University. Part of the work was done while the second and third authors visited Professor Jiayu Li in University of Science and Technology of China, they would like to thank Professor Li and the university for their hospitality and good working conditions.

\end{document}